\documentclass{amsart}

\usepackage{amssymb}
\usepackage[all]{xy}
\usepackage{hyperref}

\usepackage{enumitem}   %from SK, fixes \enumerate

\setlist[enumerate]{itemsep=.2em,topsep=.2em,leftmargin=1.25em,itemindent=2.0em}

\newtheorem{thm}{Theorem}%[section]

\newtheorem{lem}[thm]{Lemma}

\newtheorem{prop}[thm]{Proposition}

\newtheorem{complem}[thm]{Complement}%!!!!!!!!!!!!!!!!!!!!!!
\theoremstyle{definition}
\newtheorem{defn}[thm]{Definition}

\newtheorem{say}[thm]{}
\newtheorem{exmp}[thm]{Example}

\newtheorem{rem}[thm]{Remark}          

\newtheorem*{ack}{Acknowledgments}      % \renewcommand{\theack}{} 
\newtheorem{notation}[thm]{Notation}   
  
\newtheorem{defn-thm}[thm]{Definition--Theorem}  %!!!!!!!!!!!!!!!!!!!!!!!!
\newtheorem{defn-lem}[thm]{Definition--Lemma}  %!!!!!!!!!!!!!!!!!!!!!!!!
\theoremstyle{remark}

\renewcommand{\c}[0]{{\mathbb C}}

\newcommand{\z}[0]{{\mathbb Z}}

\renewcommand{\r}[0]{{\mathbb R}}

\newcommand{\dd}[0]{{\mathbb D}}

\newcommand{\q}[0]{{\mathbb Q}}

\newcommand{\qtq}[1]{\quad\mbox{#1}\quad}

\newcommand{\aut}[0]{\operatorname{Aut}}

\newcommand{\sing}[0]{\operatorname{Sing}}    
\newcommand{\ex}[0]{\operatorname{Ex}}

\newcommand{\onto}[0]{\twoheadrightarrow}

\newcommand{\tsum}[0]{\textstyle{\sum}}

\def\into{\DOTSB\lhook\joinrel\to}

\def\loccoh#1.#2.#3.#4.{H^{#1}_{#2}(#3,#4)}

\DeclareMathAlphabet{\mathchanc}{OT1}{pzc}%
                                {m}{it}

\newcommand{\GL}{\mathrm{GL}}

\usepackage[all]{xy}\xyoption{dvips}

\begin{document}
\bibliographystyle{amsalpha}

\hfill\today

 \title[Orbifold modifications]{Orbifold modifications of\\ complex analytic spaces}
 \author{J\'anos Koll\'ar and Wenhao Ou}

 \begin{abstract}
   We show that a compact, complex analytic space $X$ has a bimeromorphic orbifold modification
   that is an isomorphism over the  locally trivial orbifold locus of $X$.
 \end{abstract}

 \maketitle

 A complex analytic space $X$ is an {\it orbifold}
 if each point  has an open neighborhood  that is  biholomorphic to an  open neighborhood of  $0\in \c^n/G$ for some finite subgroup $G\subset \GL_n(\c)$. We can always assume that $G$ contains no pseudo-reflections, and then such a $G$ is unique.
 Any orbifold $X$ is {\it locally trivial} along a dense, Zariski  open subset of $\sing (X)$, that is, locally biholomorphic to an  open subset of 
$\c^r/G\times \c^{n-r}$, where $\c^r/G$ has an isolated singularity at the origin.
 Our aim is to prove the following theorem.

\begin{thm}\label{main.thm}
  Let $X$ be a compact, complex space and $S^\circ\subset \sing(X)$ a
   Zariski  open subset such that
  $X$ is a locally trivial orbifold along $S^\circ$.  Then there is a projective, bimeromorphic modification
  $\pi_Y:Y\to X$, such that  $Y$ is an orbifold, and $\pi_Y$
  is an isomorphism over both  $S^\circ$ and the smooth locus of $X$.
\end{thm}

\begin{complem} The $\pi_Y:Y\to X$ we construct has a $\pi_Y$-ample, exceptional, anti-effective divisor, and it is functorial for biholomorphisms
  (\ref{res.intro}).
 
\end{complem}

  The case when $X$ is  projective is due to  Chenyang~Xu, without assuming local triviality of the  orbifold structure along $S^\circ$. The proof  is written up in \cite[Sec.3]{MR3996488}.
  The method appears to use  projectivity in an essential way, and the 
  complex analytic version without the  local triviality assumption remains open.

Applications of Theorem~\ref{main.thm} to Bogomolov-Gieseker--type inequalities on complex spaces are treated in \cite{ou-b-g}.
  
  The proof relies on 
  some general observations on neighborhoods of subspaces and on modifications. Both of these may be of independent interest. 
  The arguments work more generally if $X$ is countable at infinity, and
  the closure of $S^\circ$ is compact.

The following example illustrates some of the difficulties, and also the method to overcome them.

\begin{exmp}
     This singularity $X:=(xy=sz^2)$ is not $\q$-factorial; $D:=(x=s=0)$ is a divisor that is not $\q$-Cartier. Thus $X$ is not a quotient singularity.

     Blowing up the origin, in the chart
     $x_1=x/s, y_1=y/s, z_1=z/s$ the equation becomes
     $x_1y_1=sz_1^2$.
     This is  the same singularity as before, so it seems that it may not have an orbifold modification, without blowing up the $s$-axis.

However, 
     adjoining  $u=\sqrt{x}, v=\sqrt{y}, t=\sqrt{s}$ we get a   Galois extension  $\bar X\to X$ with group $G=\z/2\times \z/2$, and  equation  $\bar X=(uv=tz)$.
   This can be resolved by 
  blowing up the origin to get
   $\tilde X\to \bar X\to X$.  Note that  $\tilde X$ is smooth and the
  $G$-action on $\bar X$ lifts to  a $G$-action on $\tilde X$.

  Thus  $Y:=\tilde X/G$ has only quotient singularities,
  $\pi:Y\to X$ is projective, birational, and an
  isomorphism along the $s$-axis (except at the origin).
  \end{exmp}

\begin{notation}
  All complex spaces are assumed reduced and countable at infinity.
  Note that even if we start with a normal space $X$, some of the intermediate blow-ups in our construction can be  non-normal (but reduced).  This does not cause any problems.

  Let $X$ be a complex space and $D$  a divisor on $X$.
  Then $(X, D)$ has 
  {\it simple normal crossing} at $x\in X$ if $X$ is smooth at $x$, there are local coordinates $x_1, \dots, x_n$ such that  $D=(x_1\cdots x_m=0)$ near $x$ for some $m\leq n$, 
  and the $(x_i=0)$ lie on different irreducible components of $D$.  
This is usually  abbreviated as {\it snc.}

  $\dd$ denotes a  complex disc, that is allowed to shrink as needed.
  The punctured disc is denoted by $\dd^\circ$.

  We usually work with the Euclidean  topology on $X$, and make  clear when the Zariski topology is used.
  $X^\circ$ always denotes a dense,  Zariski open subset of a complex space $X$.
\end{notation}

\begin{say}[Outline of the proof]\label{steps.plan}
 Let $X$ be a reduced,  complex space. Write $\sing(X)=S\cup R$, where $S,R$ are both 
   unions of irreducible components of $\sing(X)$, and $S$ compact.

  The proof of Theorem~\ref{main.thm} has 5 steps. In the first 2 we make no  assumptions on the singularities of $X$ along $S$. 

  \medskip
  \ref{steps.plan}.1 (Well prepared  neighborhoods)
  Let $B\subset S$ be the set of points where $X$ is not equisingular along $S$ (\ref{equising.say}).
  If $p\in S\setminus B$, then a transversal slice of $X$ through $p$ is topologically a cone in some neighborhood of radius $\epsilon(p)$.
  We use the  {\L}ojasiewicz inequality to bound $\epsilon(p)$ from below
  as $p$ goes to  $B$.
  Once this is done, it is easy to see that, 
  after a sequence of blow-ups whose centers lie over
$B$, we may assume that $S$ is smooth, there is a divisor
$D\subset X$, a Euclidean open $S\subset N\subset X$, and
 a  real analytic retraction    $\rho:N\to  S$, whose restriction
 $\rho^\circ:(N\setminus D)\to  (S\setminus D)$  is a topological fiber bundle. (There may be different fibers over different connected components of $S\setminus D$.) 
  This is done in Section~\ref{sec.1}.

   \medskip
\ref{steps.plan}.2 (Separating $S$ and $R$)
 We construct a sequence of blow-ups
  $\pi:X^r\to X$ such that 
  all blow-up centers lie over $R$,  $\pi$ induces a log resolution  of   $X\setminus S$, and
  one can write $\ex(\pi)=E^r+F^r$ such that
  $$
  \pi(E^r)\subset S\cap  R\qtq{and}
  (X^r, E^r+F^r)\qtq{is snc along}F^r.
  $$
  Note that $S^r\cup E^r=\pi^{-1}(S)$ can be very singular.
  By contrast, $F^r$, which is the closure of $\pi^{-1}(R\setminus S)$, is snc.
  Although $S^r\cup E^r$ and  $F^r$ are not disjoint, they
  only meet along the snc locus. Thus we separated the effects of $S$ from the effects of $R$. We will obtain $Y$ as a modification $Y\to X^r$ that is an isomorphism over $X^r\setminus E^r$ 

  This is done in Section~\ref{sec.2}.
We will actually need a version that also keeps track of the divisor $D$; this is small technical issue.
\medskip

Assume from now on that $X$ has an isolated quotient singularity
$\c^r/G$ transversally along $S\setminus D$.

   \medskip
 \ref{steps.plan}.3 (Local charts, first  attempt)
  Let $(V, V\cap D)\subset S$ be a local chart
  biholomorphic to $\bigl(\dd^{n-r}_{\mathbf s}, (s_1\cdots s_m=0)\bigr)$
  for some $m\leq n-r$.  Set $V^\circ:=V\setminus  D$.

  Let $U:=\rho^{-1}(V)$ be its preimage in $N$.
  Note that
  $$
  U^\circ:=U\setminus (D\cup S)=\rho^{-1}(V^\circ)\setminus S
  $$ is a topological 
  $(\c^r/G)\setminus \{0\}$ bundle over $V^\circ$.
  Its fundamental group  is an extension
  $$
  1\to G\to \pi_1(U^\circ)\to \z^m\to 0.
  $$

  Let $H=H(U)$ be a quotient of $\pi_1(U^\circ)$ such that the
  induced map $G\to H$ is an injection.
  Let $\bar U\to  U$ be the corresponding finite, normal cover.
  $\bar U$ is clearly smooth over $U\setminus (D\cup S)$,
  but it is also smooth over $S\setminus D$ since  $G\to H$ is an injection. However,   $\bar U$ can be very singular along the preimage of $D$.

  Let us now take an $H$-equivariant resolution
  $\tilde U\to \bar U$.
Then  $Y_U:=\tilde U/H \to U$
is a projective modification that
has only quotient singularities.

It is easy to see that if the quotients $\pi_1(U^\circ)\onto H$ are chosen carefully, then
the  $Y_U$ can be patched together to get the required resolution over the neighborhood $N\supset S$.
If $S$ is a connected component of $\sing X$,
then the rest of $\sing X$  can be resolved separately, and we are done.

Otherwise $\sing X$ has an irreducible component $R_i$ that
intersects both $N$ and $X\setminus N$. Then the above procedure resolved
$R_i\cap N$, but did nothing  to $R_i\setminus N$.  There is no reason to believe that a resolution along $R_i\cap N$ can be extended to one along the whole $R_i$.

We will use  (\ref{steps.plan}.2) to `isolate' $S$ from the rest of $\sing X$.  
  
  \medskip
  \ref{steps.plan}.4 (Local resolution charts)
  Let $\pi:X^r\to X$ be as in (\ref{steps.plan}.2).
  Using the local chart $U$ from (\ref{steps.plan}.3), 
set $U^r:=\pi^{-1}(U)$ and  
let $\sigma_U:\bar U^r\to  U^r$ be the  finite, normal  cover corresponding to $H$.

  Let us now take an $H$-equivariant, functorial resolution 
  $\tilde U_F^r\to \bar U^r\setminus \sigma_U^{-1}(F^r)$.  By removing
  $\sigma_U^{-1}(F^r)$, we completely disconnect $S$ from the rest of $\sing X$.

  By construction $\tilde U_F^r/H \to U^r\setminus F^r$
is a projective modification that
has only quotient singularities.

Using that $E^r$ (which is the branch locus of $\bar U^r\to U^r$) is snc along $F^r$,
an easy argument (\ref{extend.lem}) shows that  $\tilde U_F^r/H$
naturally extends to a
projective modification  $Y_U \to U^r$
 with 
only quotient singularities; see  Section~\ref{sec.3} for details.
Composing with $\pi$ we get 
$Y_U \to U$.

  \medskip
 \ref{steps.plan}.5 (Patching)
With a suitable choice of the $H$  (\ref{char.sg.say}),
the $Y_U \to U$ and
$X^r\setminus \pi^{-1}(S)\to X\setminus S$ 
patch together to the required  $Y\to X$; see Section~\ref{sec.4}.
\end{say}

\begin{say}[Similar results] In principle the method applies to other quotients, if we know enough about resolutions of their universal covers.

  For example, assume that $S$ has codimension $\geq 3$ and transversal singularity type  $(\mbox{isolated double point})/G$.
  Steps (\ref{steps.plan}.1--2) apply. In (\ref{steps.plan}.4) we get
  $\bar U^r$ that has only double points  over $V^\circ$.
  Since the resolution methods proceed by multiplicity,
  we can choose  $\tilde U_F^r\to \bar U^r\setminus \sigma_U^{-1}(F^r)$
  such that $\tilde U_F^r$ has at most double points.

  The rest of the proof works. The end result is that
  if $X$ has only double point quotients along $S^\circ$, then there is a projective, bimeromorphic modification
  $\pi_Y:Y\to X$, such that  $Y$ has only double point quotients, and $\pi_Y$
  is an isomorphism over both  $S^\circ$ and the smooth locus of $X$.

  Terminal singularities of 3-folds are double point quotients.
  Unfortunately, it is not clear that the above $Y$ also has terminal or canonical singularities, so this claim may not be very useful.
  \end{say}

\section{Well prepared neighborhoods}\label{sec.1}

Let us start with an example illustrating  well prepared  neighborhoods as in 
(\ref{steps.plan}.1).

\begin{exmp} Consider the family of curve singularities
  $$
  p: X:=(y^2=x^5+tx^2) \to \c_{t}.
  $$
  The fiber $X_t$ has an isolated singularity at the origin for every $t$, and  $X_0$ is topologically a cone.

  $X_t$ is   locally a cone for $t\neq 0$,
  but  $X_t\cap \{(x,y): |x|\leq \epsilon\}$ is a cone only for
  $\epsilon< |t|^{1/3}$.
So there is no $\epsilon>0$ such that
  $$
  X\cap \{(x,y,t): |x|<\epsilon, 0<|t|<\epsilon\}
  $$
  is a fiber bundle over $\{t: 0<|t|<\epsilon\}$.

  If we blow up the origin, then in the chart $\bar x:=x/t,  \bar y:=y/t$ the projection becomes
  $$
  \bar p: \bar X:=(\bar y^2=t\bar x^2(1+t^2\bar x^3) \to \c_{t}.
  $$
Now  $\bar X_t\cap \{(\bar x,\bar y): |\bar x|\leq \epsilon\}$ is a cone  for
$\epsilon< |t|^{-2/3}$. Thus
$$
  \bar X\cap \{(\bar x,\bar y,t): |\bar x|<1, 0<|t|<1\}
  $$
  is a fiber bundle over the punctured unit disc $\dd^\circ_t$.
\end{exmp}

\begin{say}[Equisingularity]\label{equising.say}
Let $M$ be a compact manifold. The {\it cone}  $C(M)$ over $M$ is
obtained by collapsing  $M\times \{0\}\subset M\times [0,1]$ to a  point.

Let $0\in Y\subset \c^N$ be a reduced complex space
with an isolated singularity at $0$.
Let $S(\epsilon)\subset B(\epsilon)\subset\c^N$ denote the sphere (resp.\ ball) with radius $\epsilon$ centered at $0$.

If $Y\cap S(\eta)$ is smooth for all $0<\eta\leq \epsilon$,
then  $Y\cap B(\epsilon)$ is homeomorphic to the
cone over the {\it link,}  which is  $Y\cap S(\epsilon)$.
(The homeomorphism is real analytic away from the origin.)

Let $X$ be a reduced complex space of dimension $n$, $S\subset X$ a closed subspace  of dimension $r$,
and $s\in S$ a smooth point. We say that  $X$ is
{\it equisingular}  with link $M$ along $S$ at $s$,
if there is an open neighborhood  $s\in U\subset X$, an embedding
$U\subset \dd^N_{\mathbf x}\times \dd^{r}_{\mathbf s}$
and $\epsilon>0$ such that
\begin{enumerate}
\item  $S\cap U=\{0\}\times \dd^{r}_{\mathbf s}$,
\item $U_p\cap  S(\eta)$ is smooth for all $0<\eta\leq \epsilon$ for every $p\in \dd^{r}_{\mathbf s}$, and
  \item $M$ is diffeomorphic to $U_p\cap  S(\epsilon)$.
  \end{enumerate}
If these holds then
$U\cap  \bigl(B(\epsilon)\times \dd^{r}_{\mathbf s}\bigr)$   is diffeomorphic to
$C(M)\times \dd^{r}_{\mathbf s}$.

Whitney stratification \cite{MR192520} shows that
if $X$ is a reduced complex space, then there
is a dense, Zariski open subset  $\sing^\circ X\subset \sing X$ such that
$X$ is  equisingular  along $\sing^\circ X$,  with a link that depends on the
irreducible component of $\sing^\circ X$.
\end{say}

Our aim is to show that, after some blow-ups, there is a global neighborhood that shows equisingularity. 

\begin{prop}\label{tubular.prop}
  Let $X$ be a complex space, 
  $S\subset X$ a compact subset that is a union of irreducible components of $\sing X$, and $S^\circ\subset S$ a Zariski open subset with complement $B:=S\setminus S^\circ$.
  Assume that $X$ is equisingular along $S^\circ$.

  Then there is a sequence of blow-ups
  $\pi_1:X_1\to X$ such that the following hold.
    \begin{enumerate}
\item All blow-up centers lie over $B$, so $\pi_1$ is an isomorphism over
  $X\setminus B$. In particular, $S^\circ_1:=\pi_1^{-1}(S^\circ)\cong S^\circ$.
\item The birational transform $S_1$ of $S$ is smooth, and
  $B_1:=S_1\setminus S^\circ_1$ is an snc divisor.
\item There is  a $\pi_1$-exceptional divisor $D_1$ such that
  $D_1\cap S_1=B_1$  (as schemes).
\item There is  a Euclidean neighborhood  $S_1\subset N_1\subset X_1$
  and a real analytic retraction  $\rho_1:N_1\to S_1$ 
  such that $N_1\setminus D_1$ is a locally trivial
    fiber bundle over  $S^\circ_1$.
\end{enumerate}
\end{prop}

Proof.  First we resolve the singularities of $(S, B)$
by a sequence of blow-ups whose centers lie over $B$.
(This may create new singularities outside $S$.)
Thus  it is enough to prove Proposition~\ref{tubular.prop}
when  $S$ is smooth and
$B\subset S$ is an snc divisor.

By \cite{whi-bru, MR236418} there is a real analytic embedding
$\tau:X\into \r^N$. Then $\tau(S)\subset \r^N$ has a neighborhood that is
isomorphic to the normal bundle. Its preimage by $\tau$ gives a
 Euclidean neighborhood  $S\subset N\subset X$
  and a real analytic retraction  $\rho:N\to S$

  By assumption, for each point $p\in S^\circ$ the intersection
  $\rho^{-1}(p)\cap B(\eta)$ is a cone over $M$ for $0<\eta\ll 1$.
  Let  $\epsilon(p)$ be the largest such that this holds for all
  $\eta<\epsilon(p)$. Note that (\ref{tubular.prop}.4) holds
  if $\epsilon(p)$ is uniformly bounded from below on $ S^\circ$.
The latter may not hold in  general, and 
   we need to estimate how  $\epsilon(p)$ behaves as
   $p$ approaches the boundary $B$. This is a local question along $B$.

   For $q\in B$ choose a Euclidean  neighborhood $U\ni q$ such that 
$$
U\subset \dd^M_{\mathbf x}\times \dd^r_{\mathbf s},
\qtq{with} S\cap U=\{0\}\times \dd^r_{\mathbf s}.
$$
We may also assume  that $B=(s_1\cdots s_m=0)$ for some $m\leq r$. We check in (\ref{nbhd.lem})
   that  there is an $N>0$ such that
   $$
   \epsilon(s_1,\dots, s_r)\geq |s_1\cdots s_m|^N \qtq{in a neighborhood of $q\in S$.}
   $$
 Note that if  blow up  $B$,
 the $x_i$-coordinates are replaced by  $x'_i:=x_i/(s_1\cdots s_m)$.
 Thus $\epsilon(s_1,\dots, s_r)$ changes to $\epsilon(s_1,\dots, s_r)/|s_1\cdots s_m|$.

 Choose $N$ that works for all $q\in B$.
 After $N$ blow-ups, we get $S_1\subset X_1$ and
 now $\epsilon(p)$ is uniformly bounded from below on $ S^\circ_1\cong S^\circ$.  Thus (\ref{tubular.prop}.4) holds, and the
 exceptional divisor of the last blow-up gives the  required divisor $D_1$. \qed

\begin{lem}\label{nbhd.lem} Let  $X\subset \r^m_{\mathbf x}\times \r^r_{\mathbf s}$ be the germ of a real analytic subset that contains $\{0\}\times \r^r_{\mathbf s}$.  Let $g({\mathbf s})$ be a real analytic function.
  Assume that the fibers of the coordinate projection $\pi: X\to \r^r_{\mathbf s}$ have isolated singularities at the origin $0\in \r^m_{\mathbf x}$ over $\r^r_{\mathbf s}\setminus (g=0)$. Then there is an $N$ such that
  $$
  X\cap \bigl\{({\mathbf x}, {\mathbf s}):  |{\mathbf x}|^2\leq g({\mathbf s})^N\bigr\}
  \eqno{(\ref{nbhd.lem}.1)}
  $$
  is a fiber bundle over $\r^r_{\mathbf s}\setminus (g=0)$.
\end{lem}

Proof. Fix $p\in \r^r_{\mathbf s}\setminus (g=0)$.
The fiber $X_p$ is then locally a cone over the link
$L_p:=X_p\cap S(\eta)$ for $0<\eta\ll 1$.   If the links are smooth for all
$0<\eta\leq \epsilon$, then $X_p\cap B(\epsilon)$ is a cone over
$L_p$.

Let $C'(X, \pi)$ denote the critical set
$$
\bigl\{({\mathbf x}, {\mathbf s}): X_{\mathbf s}\cap S(\eta)\mbox{ is singular at } ({\mathbf x}, {\mathbf s})\bigr\}.
$$
It is a real analytic subset of $\r^m_{\mathbf x}\times \r^r_{\mathbf s}$.
Let $C(X, \pi)$ denote the closure of  $ C'(X, \pi)\setminus \bigl(\{0\}\times \r^r_{\mathbf s}\bigr)$.

Then  $C(X, \pi)$
 is a real semianalytic set by \cite[2.8]{MR972342},
 and  $g$ vanishes at all points of $ C(X, \pi)$ where
 $\tsum |x_i|^2$ vanishes.
 By the {\L}ojasiewicz inequality
 (see \cite{loj} or
\cite[6.4]{MR972342})
 there is an $N$ such that
 $$
 \tsum |x_i|^2|_{C(X, \pi)}\geq g({\mathbf s})^N|_{C(X, \pi)}. \qed
  \eqno{(\ref{nbhd.lem}.2)}
  $$

 \section{Constrained resolutions}\label{sec.2}

 We discuss how much can be resolved, if we are not allowed to blow up over some subset $S^\circ$.

 \begin{defn}\label{res.intro}
   Let $X$ be a complex space.
   A {\it modification} is a proper, bimeromorphic morphism $\pi:Y\to X$.
   It is a {\it resolution} if $Y$ is smooth.

   A modification is 
   {\it functorial for biholomorphisms} if the following holds.
   Let $U_1, U_2\subset X$ be open, and
  $U_1\cong U_2$  a biholomorphism. Then it induces a biholomorphism
   $\pi^{-1}(U_1)\cong \pi^{-1}(U_2)$.

   Note that such a modification is necessarily an isomorphism over the smooth locus of $X$.

   The existence of a functorial resolution goes back to \cite{hir};
   see \cite{k-res}
   for an introduction and proofs.
   \end{defn}

 \begin{lem}\label{constr.prop.1}
   Let $X$ be a complex space,  $D\subset X$ a divisor, and
   $\sing(X, D)$  the closed subset where $X$ or $D$ is singular.
   Write $\sing(X, D)=S\cup R$, where $S,R$ are both 
   unions of irreducible components of $\sing(X, D)$.
   Set $\bar R:=R\cup (S\cap D)$.
   Assume that $S\cap D\subset S$ is nowhere dense.

  Then there is a sequence of blow-ups
  $\pi:X^r\to X$ such that 
  \begin{enumerate}
\item all blow-up centers lie over $\bar R$ (so $\pi$ is an isomorphism over
  $X\setminus \bar R$), and
\item  $\pi$ induces a log resolution  of   $(X\setminus S, D\setminus S)$. 
  \end{enumerate}
  Moreover, if $S$ is smooth and $B\subset S$ is an snc divisor, then
  we can choose the sequence of blow-ups such that
    \begin{enumerate}\setcounter{enumi}{2}
    \item $S^r:=\pi^{-1}_*(S)$ is smooth, and $B^r:=S^r\cap \pi^{-1}(B)$ is an snc divisor.
        \end{enumerate}
\end{lem}

 Proof.  Set $(X_0, D_0):=(X, D)$ and
 let $\{\tau_i: X_{i+1}:=B_{Z_i}X_i\to X_i\}$ be a resolution  sequence for
 $(X, D)$.   Let $\Pi_i:X_i\to X_0$ denote the composites.

 Set $X^*:=X\setminus S$.   By restriction we get a
 resolution  sequence for  $(X^*, D^*)$.
 We extend it to a blow-up sequence of $X$ as follows.

 Assume that we already have  $\Pi'_i:X'_i\to X'_0=X$
 which is isomorphic to  $\Pi^*_i:X^*_i\to X^*_0$ over $X^*$.

 If $Z^*_i:=Z_i\cap X^*_i$ is nonempty, let $Z'_i\subset X'_i$ be its closure.
 Set $X'_{i+1}:=B_{Z'_i}X'_i$.
 If  $Z^*_i$ is empty, we set $X'_{i+1}:=X'_i$.
This shows (\ref{constr.prop.1}.1--2).

For (\ref{constr.prop.1}.3) we need to modify the above sequence as follows.
First blow up over $S\cap R$ to achieve
that $S$ is smooth and $B$ is an snc divisor in $S$.
Then we run the above procedure, with aditional blow-ups.

Note that even if the $Z_i$ are smooth, the $Z'_i$ can be very singular, and the smoothness of $S$ is not preserved.
To remedy this, assume that we already have  $\Pi'_i:X'_i\to X'_0$ and 
 $S'_i\subset X'_i$ (the birational transform of $S$) is smooth. Then we can perform first an  auxiliary sequence of blow-ups
over $S'_i\cap Z'_i$  to get
$X''_i\to X'_i$ such that $(S''_i, B''_i)$ is snc, and
$Z''_i$ meets $(S''_i, B''_i)$ transversally.

We then let $X'_{i+1}$  be the blow-up of $Z''_i\subset X''_i$. \qed

\medskip

 Note that $\ex(\pi)$  is an snc divisor over $X\setminus S$,
but we have no control over its singularities that lie over $S$.
We need a stronger version that takes care of this problem.

 \begin{prop}\label{constr.prop.2}
    Using the notation and assumptions of (\ref{constr.prop.1}), 
 there is a sequence of blow-ups
  $\pi:X^r\to X$ such that 
  all blow-up centers lie over $\bar R$,  $\pi$ induces a log resolution  of   $(X\setminus S, D\setminus S)$, and
  one can write $\ex(\pi)=E^r+F^r$ such that
  \begin{enumerate}
  \item  $\pi(E^r)\subset S\cap \bar R$, and
    \item  $(X^r, E^r+F^r+D^r)$ is snc along  $F^r+D^r$.
\end{enumerate}
\end{prop}

 Proof. Let
 $\tau: X_1\to X$  be the blow-up of the (ususally non-reduced)  $S\cap \bar R$. Then   
   $S_1:=\tau^{-1}_*(S)$ is disjoint from
 $\tau^{-1}_*(R)$ and $D_1:=\tau^{-1}_*(D)$.
 Then we apply (\ref{constr.prop.1}) to
 $X_1, S_1$ and $D_1$. \qed

 \begin{rem} With $X$ and $S$ as in Theorem~\ref{main.thm},
   we will construct $Y$ as a modification
   $Y\to X^r$ that is an isomorphism over
   $X^r\setminus E^r$. Thus we no longer need to pay attention
   to the singularities in $R$.
   \end{rem}

 \section{Resolution of local charts}\label{sec.3}

  \begin{say}[Local construction of the resolution]\label{char.sg.say}
Following the outcome of Proposition~\ref{tubular.prop},
  let $X$ be a complex space, 
  $S\subset X$ a compact subset that is a union of irreducible components of $\sing X$,   $D\subset X$ a divisor and
  a Euclidean neighborhood  $S\subset N\subset X$. We also have 
  a real analytic retraction  $\rho:N\to S$,
  whose restriction
  $$
  \rho^\circ:N^\circ:=(N\setminus (D\cup S))\to S^\circ:=(S\setminus D)
  $$
  is  a fiber bundle, and the  fiber has finite fundamental group $G$.
  We also assume that  $(S, S\cap D)$ is an snc pair.

  Let $\pi:X^r\to X$ be as in (\ref{constr.prop.2}). 

  Choose a local chart  $V\cong\dd^{a+c}$ on $S$ such that  $V^\circ=\dd^a\times (\dd^\circ)^c$.   Let $U\subset N$ be its preimage and set
  $U^\circ:=U\setminus (D\cup S)$. 
  We get a fiber bundle
$U^\circ\to V^\circ$.
  Since $\pi_2(V^\circ)=0$, there is a  group extension
  $$
  1\to G\to \pi_1(U^\circ)\to \pi_1(V^\circ)\cong\z^c\to 0.
  \eqno{(\ref{char.sg.say}.1)}
  $$
  Let $\pi_1(U^\circ)\onto H$ be a quotient such that
  $G\to H$ is an injection, and $\bar U\to U$ the corresponding finite cover that is 
   \'etale over $U^\circ$, but usually ramified over $U\cap (D\cup S)$.

  Let $U^r:=\pi^{-1}(U)\subset X^r$ be the preimage of $U$ and
  $\sigma_U:\bar U^r\to U^r$ the corresponding cover.
We remove the preimage $\sigma_U^{-1}(F^r\cup D^r)$ from
$\bar U^r$.  Doing this we completely disconnected $S$ from the rest of $\sing X$.

Assume form now on that   $X$ has quotient singularities along $S^\circ$. Then
  $\bar U$  and   $\bar U^r$ are smooth along the preimage of $S^\circ$.

  Next take an $H$-equivariant, functorial resolution 
  $\tilde U_F^r\to \bar U^r\setminus \sigma_U^{-1}(F^r\cup D^r)$.  

  By construction $\tilde U_F^r/H \to U^r\setminus (F^r\cup D^r)$
is a projective modification that
has only quotient singularities.

Note that $\bar U^r\to U^r$ is ramified only along
$E^r\cup F^r\cup D^r$. Although $E^r$ is very singular,
the ramification locus is snc along the intersection
$E^r\cap (F^r\cup D^r)$. Thus we can apply 
(\ref{extend.lem}) to the snc locus to extend  $\tilde U_F^r/H$
 to a
projective modification  $Y_U \to U^r$
 with 
only quotient singularities. Composing with $\pi$ we get 
$Y_U \to U$.
\end{say}

\begin{lem}\label{extend.lem}
  Let  $(Z, E+F)$ be an snc pair, and
   $\tau_H: Z_H\to Z$ a finite Galois cover with group $H$, ramified along $E+F$.  Let $\tilde Z_H\to Z_H\setminus \tau_H^{-1}(F)$ be a functorial modification.  Then $\tilde Z_H/H\to (Z\setminus F)$ has a unique extension to a  functorial modification
   $Y\to Z$ of $(Z, E)$. In particular, $Y\to Z$  is an isomorphism over $Z\setminus E$.
 \end{lem}

 Proof. The question is local on $Z$, so we may assume that
 $Z=\dd^n$,  $E=(x_1\cdots x_r=0)$ and $F=(x_{r+1}\cdots x_n=0)$.
 Let $K\subset H$ be the image of $\pi_1(Z\setminus E)$, and
 $\tau_K:Z_K\to Z$ the corresponding cover.
 The induced map $Z_H\to \bigl(Z_K\setminus \tau_K^{-1}(F)\bigr)$ is a finite, \'etale cover.  Thus $\tilde Z_H$ decends to
 a functorial modification $\tilde Z_{F,K}$ of $Z_K\setminus \tau_K^{-1}(F)$.
 The product structure on $Z_K$ then extends this to
a functorial modification  $\tilde Z_K\to Z_K$. 
Set $Y:=\tilde Z_K/K$. \qed

\section{Patching  the local charts}\label{sec.4}

\begin{say}[Patching]\label{patch.say}
Continuing with the notation of
(\ref{char.sg.say}), choose local charts $\{V_i:i\in I\}$ that cover  an irreducible component of
$S$. We have  $U^\circ_i\subset N$ and
 fiber bundles
$U_i^\circ\to V_i^\circ$. As in (\ref{char.sg.say}.1), 
  there are  group extensions
  $$
  1\to G\to \pi_1(U_i^\circ)\to \pi_1(V_i^\circ)\cong\z^{c_i}\to 0.
  %%\eqno{(\ref{char.sg.say}.1)}
  $$
  For each $i$ we choose a quotient
  $\pi_1(U_i^\circ)\onto H_i$   such that the induced homomorphisms
  $G\to H_i$ are injections.
  We get the ramified covers  $\bar U_i\to U_i$.

  Since the construction of the resolutions $U_i^r$, the covers
  $(\bar U_i)_F^r$ and the resolutions $(\tilde U_i)_F^r$ are all
  functorial for \'etale morphisms, the $Y_{U_i}\to U_i$
  patch together if the fiber product maps
  $$
  \bar U_i \leftarrow \bar U_i\times_S \bar U_j \rightarrow \bar U_j
  \eqno{(\ref{patch.say}.1)}
  $$
  are \'etale.
  We do not understand the singularities of the $\bar U_i$, but
  this holds if 
  the fiber product maps
  $$
  \bar U_i^\circ \leftarrow \bar U_i^\circ\times_S  \bar U_j^\circ \rightarrow  \bar U_j^\circ
  \eqno{(\ref{patch.say}.2)}
  $$
  are open embeddings on each connected component.
  
  Let $\bar V_i\to V_i$ be the ramified cover corresponding to $H_i/G$.
  Since the  
  $\bar U_i^\circ\to \bar V_i^\circ$ are  fiber bundles with the same fibers, the projections in 
  (\ref{patch.say}.2) are open embeddings  iff the
  $$
  \bar V_i^\circ \leftarrow \bar V_i^\circ\times_S  \bar V_j^\circ \rightarrow  \bar V_j^\circ
  \eqno{(\ref{patch.say}.3)}
  $$
  are open embeddings on each connected component.

  If the $\bar V_i\to V_i$  ramify along each irreducible component
  of $V_i\cap D$ with the same ramification index $m$, then
  the projections in 
  (\ref{patch.say}.3) are \'etale.
  In (\ref{char.split.say}) we choose $m$ and the $H_i$ 
  such that the covering group of $\bar V_i^\circ\to V_i^\circ$ is
  $H_1(V_i^\circ, \z)/m H_1(V_i^\circ, \z)$ for every $i$.
  Such a choice of the $H_i$ is described next.

  If the intersections of 
  $V_i\cap V_j$  with the strata of $S\cap D$ are simply connected,
  then $H_1(V_i^\circ\cap V_j^\circ, \z)\to H_1(V_i^\circ, \z)$
  is a direct summand.  Then the restriction of
  $\bar V_i^\circ\to V_i^\circ$ to $V_i^\circ\cap V_j^\circ$
  splits into a disjoint union of the
  $H_1(V_i^\circ\cap V_j^\circ, \z)/m H_1(V_i^\circ\cap V_j^\circ, \z)$ covers of 
$V_i^\circ\cap V_j^\circ$. \qed

\end{say}

  \begin{say}[Characteristic splitting]\label{char.split.say} Let
    $$
  1\to G\to\Gamma\to \z^{c}\to 0
  %%\eqno{(\ref{char.sg.say}.1)}
  $$
  be a group extension.
  Set $m_1:=|\aut(G)|$. Then the  image of the centralizer of $G$  in
$\z^c$
  contains $(m_1\z)^c$. This gives a subgroup 
  $\Gamma_1\subset \Gamma$ which is a central extension
  $$
  1\to Z(G)\to \Gamma_1\to (m_1\z)^c\to 0.
  $$
  Let $m_2:=|Z(G)|$. Then the $m_2$-th power of any element of $\Gamma_1$  lies in its center. So there is an abelian subgroup
  $\Gamma_2\subset\Gamma_1$ and an extension
  $$
  1\to Z(G)\to \Gamma_2\to (m_1m_2\z)^c\to 0.
  $$
    Set $m=m_1m_2^2$. The $m_2$-th powers of elements of
  $\Gamma_2$ generate a central,  characteristic subgroup  $\Gamma_3\subset\Gamma_2\subset  \Gamma$,
  which maps isomorphically onto
  $(m\z)^c\subset \z^c$. The quotient is an extension
  $$
  1\to G\to H\to (\z/m)^c\to 0.
  $$
  \end{say}

    \begin{ack}
      We are grateful to Omprokash~Das for helpful conversations.
      Partial  financial support  to JK   was provided  by the Simons Foundation   (grant number SFI-MPS-MOV-00006719-02).
Partial  financial support  to WO   was provided  by
      the National Key R{\&}D Program of China (No.2021YFA1002300).
    \end{ack}

 %%\bibliography{../refs-main/refs}

    \def\cprime{$'$} \def\cprime{$'$} \def\cprime{$'$} \def\cprime{$'$}
  \def\cprime{$'$} \def\dbar{\leavevmode\hbox to 0pt{\hskip.2ex
  \accent"16\hss}d} \def\cprime{$'$} \def\cprime{$'$}
  \def\polhk#1{\setbox0=\hbox{#1}{\ooalign{\hidewidth
  \lower1.5ex\hbox{`}\hidewidth\crcr\unhbox0}}} \def\cprime{$'$}
  \def\cprime{$'$} \def\cprime{$'$} \def\cprime{$'$}
  \def\polhk#1{\setbox0=\hbox{#1}{\ooalign{\hidewidth
  \lower1.5ex\hbox{`}\hidewidth\crcr\unhbox0}}} \def\cdprime{$''$}
  \def\cprime{$'$} \def\cprime{$'$} \def\cprime{$'$} \def\cprime{$'$}
\providecommand{\bysame}{\leavevmode\hbox to3em{\hrulefill}\thinspace}
\providecommand{\MR}{\relax\ifhmode\unskip\space\fi MR }
% \MRhref is called by the amsart/book/proc definition of \MR.
\providecommand{\MRhref}[2]{%
  \href{http://www.ams.org/mathscinet-getitem?mr=#1}{#2}
}
\providecommand{\href}[2]{#2}

\bigskip

  Princeton University, Princeton NJ 08544-1000, 

  \email{kollar@math.princeton.edu}

  Institute of Mathematics, Academy of Mathematics and Systems Science,

  Chinese
Academy of Sciences, Beijing, 100190, China

\email{wenhaoou@amss.ac.cn}
 
\end{document}